\documentclass[12pt,  reqno]{amsart}

\newtheorem{theorem}{Theorem}
\theoremstyle{plain}

\newtheorem{corollary}[theorem]{Corollary}

\newtheorem{remark}{Remark}

\numberwithin{equation}{section}

\begin{document}
\title[] {Some identities of the generalized twisted Bernoulli numbers and polynomials \\of higher order}
\author[] {Young-Hee Kim}
\address{Young-Hee Kim. Division of General Education-Mathematics\\
Kwangwoon University, Seoul 139-701, Republic of Korea  \\}
\email{yhkim@kw.ac.kr}
\author[]{Byungje Lee}
\address{Byungje Lee. Department of Wireless Communications Engineering \\
Kwangwoon University, Seoul 139-701, Republic of Korea  \\}
\email{bjlee@kw.ac.kr}
\author[]{Seog-Hoon Rim}
\address{Seog-Hoon Rim. Department of Mathematics Education\\
Kyungpook National University, Taegu 702-701, Republic of Korea
\\} \email{shrim@knu.ac.kr}
\author[]{Taekyun Kim}
\address{Taekyun Kim. Division of General Education-Mathematics \\
Kwangwoon University, Seoul 139-701, Republic of Korea  \\}
\email{tkkim@kw.ac.kr}

\maketitle

{\footnotesize {\bf Abstract} \hspace{1mm} {The purpose of this
paper is to derive some identities of the higher order generalized
twisted Bernoulli numbers and polynomials attached to $\chi$ from
the properties of the $p$-adic invariant integral.  We give some
interesting identities for the power sums and the generalized
twisted Bernoulli numbers and polynomials of higher order using the
symmetric properties of the $p$-adic invariant integral. }



\section{Introduction and preliminaries}
Let $p$ be a fixed prime number. Throughout this paper, the symbol
$\Bbb Z$, $\Bbb Z_p$, $\Bbb Q_p$, and $\Bbb C_p$ denote the ring of
rational integers, the ring of $p$-adic integers, the field of
$p$-adic rational numbers, and the completion of algebraic closure
of $\Bbb Q_p$, respectively. Let $\mathbb{N}$ be the set of natural
numbers and $\Bbb Z_+ = \Bbb N \cup \{ 0 \}$ . Let $\nu_p $ be the
normalized exponential valuation of $\Bbb C_p$ with $|p|_p=
p^{-\nu_{p}(p)}=p^{-1}$.

Let $UD(\Bbb Z_p)$ be the space of uniformly differentiable function
on $\Bbb Z_p$. For $f \in UD(\Bbb Z_p)$, the $p$-adic invariant
integral on $\Bbb Z_p$ is defined as
\begin{eqnarray} I(f)=\int_{\Bbb Z_p
} f(x) dx= \lim_{N \rightarrow \infty} \frac{1}{p^N} \sum_{x=0}^{p^N
-1} f(x).
\end{eqnarray}
(see [4-5]). From $(1.1)$, we note that
\begin{eqnarray} I(f_1)=I(f)+f'(0),
\end{eqnarray} where $ f'(0)=\frac{df(x)}{dx} |_{x=0}$ and $f_1(x)=f(x+1)$.
For $n \in \mathbb{N}$, let $f_n(x)=f(x+n)$. Then we can derive the
following equation from (1.2).
\begin{eqnarray} I(f_n)=I(f)+ \sum_{i=0}^{n-1}f'(i), \quad (\text{see [4-5]}). \label{1.3}
\end{eqnarray}

Let $d$ be a fixed positive integer. For $n \in \mathbb{N}$, let
\begin{eqnarray*}
X&=& X_d = \lim_{\overleftarrow{N} } \Bbb Z/ dp^N \Bbb Z ,\ X_1 =
\Bbb Z_p , \\   X^\ast &=&\underset{{0<a<d p}\atop {(a,p)=1}} {\cup}
(a+ dp \,\Bbb Z_p ), \\  a+d p^N \Bbb Z_p &=& \{ x\in X | \,  x
\equiv a \pmod{dp^N}\},
\end{eqnarray*}
where $a\in \Bbb Z$ lies in $0\leq a < d p^N$. It is easy to see
that
\begin{eqnarray}\int_X f(x)dx = \int_{\Bbb Z_p} f(x) dx, \quad
\text{for} \quad f \in UD(\Bbb Z_p). \label{1.4} \end{eqnarray}

The ordinary Bernoulli polynomials $B_n (x)$ are defined as
\begin{eqnarray*} \frac{t}{ e^t -1 }e^{xt} =\sum_{n=0}^{\infty} B_n (x)
\frac{t^n}{n!},  \end{eqnarray*} and the Bernoulli numbers $B_n$ are
defined as $B_n=B_n(0)$ (see [1-19]).

For $n \in \Bbb N$, let $T_p$ be the p-adic locally constant space
defined by
$$T_p =\underset{n \ge 1} {\cup} \Bbb C_{p^n}= \underset{n \to
\infty} {\lim} \Bbb C_{p^n},$$ where $\Bbb C_{p^n}= \{ \omega |
\omega^{p^n}=1 \}$ is the cyclic group of order $p^n$. It is well
known that the twisted Bernoulli polynomials are defined as
\begin{eqnarray*} \frac{t}{\xi e^t -1 }e^{xt} =\sum_{n=0}^{\infty} B_{n, \xi}(x)
\frac{t^n}{n!}, \quad \xi \in T_p ,
\end{eqnarray*}
and the twisted Bernoulli numbers $B_{n, \xi}$ are defined as $B_{n,
\xi}=B_{n, \xi}(0)$ (see [14-18]).

Let $\chi$ be the Dirichlet's character with conductor $d \in
\mathbb{N}$. Then we have
\begin{eqnarray}
\int_{X}\chi(x) \xi^x e^{xt}dx=
\frac{t\underset{a=0}{\overset{d-1}{\sum}}\chi(a)\xi^a e^{at}}{\xi^d
e^{dt}-1} . \label{1.5}
\end{eqnarray}

It is known that the generalized twisted Bernoulli numbers attached
to $\chi$, $B_{n,\chi, \xi}$, are defined as
\begin{eqnarray}
\frac{t\underset{a=0}{\overset{d-1}{\sum}}\chi(a)\,\xi^a e^{at}
}{\xi^d e^{dt} -1}= \underset{n=0}{\overset{\infty}{\sum}}
B_{n,\chi, \xi} \frac{t^n}{n!}, \quad \xi \in T_p. \label{1.6}
\end{eqnarray}
The generalized twisted Bernoulli polynomials attached to $\chi$,
$B_{n,\chi, \xi} (x)$, are defined as
\begin{eqnarray}
\frac{t\underset{a=0}{\overset{d-1}{\sum}}\chi(a)\,\xi^a e^{at}
}{\xi^d e^{dt} -1}e^{xt}= \underset{n=0}{\overset{\infty}{\sum}}
B_{n,\chi, \xi} (x)\frac{t^n}{n!}, \quad \xi \in T_p, \label{1.7}
\end{eqnarray}
(see [13], [16]). From ($\ref{1.5}$), ($\ref{1.6}$) and
($\ref{1.7}$), we derive that
\begin{eqnarray}
\int_{X}\chi(x) \xi^x x^ndx= B_{n,\chi, \xi} \quad \text{and} \quad
\int_{X}\chi(y) \xi^y (x+y)^n dy= B_{n,\chi, \xi}(x).
\label{1.8}\end{eqnarray}

By ($\ref{1.3}$) and ($\ref{1.4}$), it is easy to see that for $n
\in \mathbb{N}$,
\begin{eqnarray}
\int_{X}f(x+n)dx=\int_{X}f(x)dx+\sum_{i=0}^{n-1}f\,'(i), \label{1.9}
\end{eqnarray}
where $f\,'(i)=\frac{df(x)}{dx}|_{x=i}$. From ($\ref{1.9}$), it
follows that

\begin{eqnarray}
&&\frac{1}{t}(\int_{X}\chi(x)\,\xi^{nd+x}e^{(nd+x)t}dx-\int_{X}\chi(x)\,\xi^{x}e^{xt}dx)
\label{1.10}\\&& = \frac{nd\int_{X}\chi(x)\,\xi^x
e^{xt}dx}{\int_{X}\xi^{ndx}e^{ndxt}
dx}=\frac{\xi^{nd}e^{ndt}-1}{\xi^d
e^{dt}-1}(\sum_{i=0}^{d-1}\chi(i)\,\xi^i e^{it})
=\sum_{k=0}^{\infty}(\sum_{l=0}^{nd-1}\chi(l)\, \xi^l l^{k})
\frac{t^k}{k!}. \quad \notag
\end{eqnarray}

For $k \in \Bbb Z_+$, let us define the $p$-adic functional $T_{k,
\chi, \xi}(n)$ as follows:
\begin{eqnarray}
T_{k, \chi, \xi}(n)=\sum_{l=0}^{n}\chi(l)\xi^l l^{k}. \label{1.11}
\end{eqnarray}
Let $k,n,d \in \mathbb{N}$. By $(\ref{1.10})$ and $(\ref{1.11})$, we
see that
\begin{eqnarray}
\int_{X} \chi(x)\xi^{nd+x}(nd+x)^k dx-\int_{X}\chi(x) \xi^x x^k dx=k
\,T_{k-1, \chi, \xi}(nd-1).\label{1.12}
\end{eqnarray}
From $(\ref{1.8})$ and $(\ref{1.12})$, we have that
\begin{eqnarray}
\frac{\xi^{nd} B_{k,\chi, \xi}(nd)-B_{k,\chi, \xi}}{k}=T_{k-1, \chi,
\xi}(nd-1). \label{1.13}
\end{eqnarray}

For $w_1, w_2, d \in \mathbb{N}$, we note that
\begin{eqnarray}
& &\frac{d\int_{X}\int_{X}\chi(x_1)\chi(x_2)\,\xi^{w_1 x_1+w_2 x_2}
e^{(w_1 x_1+w_2 x_2)t}
 dx_1
dx_2 }{\int_{X} \xi^{dw_1 w_2 x }e^{dw_1 w_2 xt} x}\label{1.14}\\
& &\quad =\frac{t(\xi^{dw_1 w_2}e^{dw_1 w_2 t} -1)}{(\xi^{w_1
d}e^{w_1d t}-1)(\xi^{w_2 d}e^{w_2 d t}-1)}
(\sum_{a=0}^{d-1}\chi(a)\xi^{w_1 a}e^{w_1 at}
)(\sum_{b=0}^{d-1}\chi(b)\xi^{w_2 b}e^{w_2 bt} ). \quad \notag
\end{eqnarray}
In the next section, we will consider the extension of
$(\ref{1.14})$ related to the generalized twisted Bernoulli numbers
and polynomials of higher order attached to $\chi$ .

The generalized twisted Bernoulli polynomials of order $k$ attached
to $\chi$, $B_{n,\chi, \xi}^{(k)}(x)$, are defined as
\begin{eqnarray}
\left( \frac{t\underset{a=0}{\overset{d-1}{\sum}}\chi(a)\,\xi^a
e^{at} }{\xi^d e^{dt} -1} \right)^k e^{xt}=
\underset{n=0}{\overset{\infty}{\sum}} B_{n,\chi, \xi}^{(k)}
(x)\frac{t^n}{n!}, \quad \xi \in T_p, \label{1.15}
\end{eqnarray}
and $B_{n,\chi, \xi}^{(k)}=B_{n,\chi, \xi}^{(k)}(0)$ are called the
generalized twisted Bernoulli numbers of order $k$ attached to
$\chi$. When $k=1$, the polynomials and numbers are called the
generalized twisted Bernoulli polynomials and numbers attached to
$\chi$, respectively (see [12]).

The authors of this paper have studied various identities for the
Bernoulli and the Euler polynomials by the symmetric properties of
the $p$-adic invariant integrals (see [6-8], [10]). T. Kim [6]
established interesting identities by the symmetric properties of
the $p$-adic invariant integrals and some relationships between the
power sums and the Bernoulli polynomials. In [8], Kim et al. gave
some identities of symmetry for the generalized Bernoulli
polynomials. The twisted Bernoulli polynomials and numbers are very
important in several field of mathematics and physics, and so have
been studied by many authors (cf. [9-18]). Recently, Kim-Hwang [10]
obtained some relations between the power sum polynomials and
twisted Bernoulli polynomials.


In this paper, we extend our results to the generalized twisted
Bernoulli numbers and polynomials of higher order attached to
$\chi$. The purpose of this paper is to derive some identities of
the higher order generalized twisted Bernoulli numbers and
polynomials attached to $\chi$ from the properties of the $p$-adic
invariant integral. In Section 2, we give interesting identities for
the power sums and the generalized twisted Bernoulli numbers and
polynomials of higher order using the symmetric properties for the
$p$-adic invariant integral.

\medskip

\section{Some identities of the generalized twisted Bernoulli numbers and polynomials of higher order}

Let $w_1, w_2, d \in \mathbb{N}$. For $\xi \in T_p$, we set

\begin{eqnarray}
& &Y(m, \chi, \xi | w_1, w_2) \notag \\
&&\quad=\left(\frac{d\int_{X^m}(\underset{i=1}{\overset{m}{\prod}}\chi(x_i))\xi^{(\underset{i=1}{\overset{m}{\sum}}
x_i)w_1 } e^{(\underset{i=1}{\overset{m}{\sum}}x_i+w_2x)w_1t}dx_1
\cdots dx_m} {\int_{X}\xi^{dw_1 w_2 x}e^{dw_1w_2xt}dx}\right)\label{2.1}\\
& & \qquad\times
\left(\int_{X^m}(\underset{i=1}{\overset{m}{\prod}}\chi(x_i))\xi^{(\underset{i=1}{\overset{m}{\sum}}
x_i)w_2 }e^{(\underset{i=1}{\overset{m}{\sum}}x_i+w_1y)w_2t} dx_1
\cdots dx_m \right), \notag
\end{eqnarray}
where
\begin{eqnarray*}
\int_{X^m}f(x_1, \cdots, x_m)dx_1 \cdots dx_m= \underbrace{\int_{X}
\cdots \int_{X}}_{m-\text{times}}f(x_1, \cdots, x_m)dx_1 \cdots
dx_m.
\end{eqnarray*}

In $(\ref{2.1})$, we note that $Y(m, \chi, \xi; w_1, w_2)$ is
symmetric in $w_1, w_2$. From $(\ref{2.1})$, we derive that
\begin{eqnarray}
& &Y(m, \chi, \xi | w_1, w_2) \notag \\
&&=\left(\int_{X^m}(\underset{i=1}{\overset{m}{\prod}}\chi(x_i))\xi^{(\underset{i=1}{\overset{m}{\sum}}
x_i)w_1 }e^{(\underset{i=1}{\overset{m}{\sum}}x_i)w_1t}dx_1 \cdots
 dx_m\right)e^{w_1w_2xt}\qquad \label{2.2}\\
& & \quad \times
 \left(\frac{d\int_{X}\chi(x_m)\xi^{w_2 x_m}e^{w_2x_{m}t}dx_m}{\int_{X}\xi^{dw_1 w_2
 x}e^{dw_1w_2xt}dx}\right) \notag\\
& & \quad \times
 \left(\int_{X^{m-1}}(\underset{i=1}{\overset{m-1}{\prod}}\chi(x_i))\xi^{(\underset{i=1}{\overset{m-1}{\sum}}
x_i)w_2 }e^{(\underset{i=1}{\overset{m-1}{\sum}}x_i)w_2t}dx_1
 \cdots  dx_{m-1}\right)e^{w_1w_2yt}. \notag
\end{eqnarray}
From $(\ref{1.10})$ and $(\ref{1.11})$, it follows that
\begin{eqnarray}
\frac{dw_1\int_{X}\chi(x) \xi^x e^{xt}dx}{\int_{X}
\xi^{dw_1x}e^{dw_1xt}dx}=\sum_{i=0}^{w_1d-1}\chi(i)\xi^i e^{it}
=\sum_{k=0}^{\infty}T_{k, \chi, \xi}(w_1d-1)\frac{t^k}{k!}.
\label{2.3}
\end{eqnarray}
By $(\ref{1.15})$, we also see that
\begin{eqnarray}
&&e^{w_1w_2xt}\left(\int_{X^m}(\underset{i=1}{\overset{m}{\prod}}\chi(x_i))\xi^{(\underset{i=1}{\overset{m}{\sum}}
x_i)w_1 }e^{(\underset{i=1}{\overset{m}{\sum}}x_i)w_1t}dx_1
 \cdots  dx_m \right)\label{2.4}\\ & & \quad
 = \left(\frac{w_1t}{\xi^{dw_1}e^{dw_1t}-1}\sum_{a=0}^{d-1}\chi(a)\xi^{w_1 a}e^{aw_1t}\right)^m e^{w_1w_2xt}=
 \underset{n=0}{\overset{\infty}{\sum}} B_{n,\chi, \xi^{w_1}}^{(m)} \notag
(w_{2}x)\frac{w_1^nt^n}{n!}.
\end{eqnarray}
By $(\ref{2.2})$, $(\ref{2.3})$ and $(\ref{2.4})$, we have that
\begin{eqnarray}
& &Y(m, \chi, \xi | w_1, w_2) \label{2.5} \\
& &=\left(\underset{l=0}{\overset{\infty}{\sum}} B_{l,\chi,
\xi^{w_1}}^{(m)}
(w_{2}x)\frac{w_1^lt^l}{l!}\right)\left(\frac{1}{w_1}\underset{k=0}{\overset{\infty}{\sum}}T_{k,
\chi,
\xi^{w_2}}(w_1d-1)\frac{w_2^kt^k}{k!}\right)\left(\underset{i=0}{\overset{\infty}{\sum}}
B_{i,\chi,
\xi^{w_2}}^{(m-1)}(w_{1}y)\frac{w_2^it^i}{i!}\right)\notag
\\
&&=\underset{n=0}{\overset{\infty}{\sum}}\
\left(\underset{j=0}{\overset{n}{\sum}}\binom{n}{j}w_2^{j}
w_1^{n-j-1}B_{n-j, \, \chi,
\xi^{w_1}}^{(m)}(w_2x)\underset{k=0}{\overset{j}{\sum}}\binom{j}{k}T_{k,
\chi, \xi^{w_2}}(w_1d-1)B_{j-k,\chi,
\xi^{w_2}}^{(m-1)}(w_{1}y)\right)\frac{t^n}{n!}.\notag
\end{eqnarray}
From the symmetry of $Y(m, \chi, \xi | w_1, w_2)$ in $w_1$ and
$w_2$, we see that
\begin{eqnarray}
& &Y(m, \chi, \xi | w_1, w_2) \label{2.6}\\
&&=\underset{n=0}{\overset{\infty}{\sum}}\left(\underset{j=0}{\overset{n}{\sum}}\binom{n}{j}w_1^{j}
w_2^{n-j-1}B_{n-j, \chi,
\xi^{w_2}}^{(m)}(w_1x)\underset{k=0}{\overset{j}{\sum}}\binom{j}{k}T_{k,
\chi, \xi^{w_1} }(w_2d-1)B_{j-k, \chi,
\xi^{w_2}}^{(m-1)}(w_2y)\right)\frac{t^n}{n!}.\notag
\end{eqnarray}

Comparing the coefficients on the both sides of $(\ref{2.5})$ and
$(\ref{2.6})$, we obtain an identity for the generalized twisted
Bernoulli polynomials of higher order as follows.
\begin{theorem}
Let $d, w_1, w_2 \in \mathbb{N}$. For $n \in \mathbb{Z}_+$ and $m
\in \mathbb{N}$, we have
\begin{eqnarray*}
& &\underset{j=0}{\overset{n}{\sum}}\binom{n}{j}w_2^{j}
w_1^{n-j-1}B_{n-j, \, \chi,
\xi^{w_1}}^{(m)}(w_2x)\underset{k=0}{\overset{j}{\sum}}\binom{j}{k}T_{k,
\chi, \xi^{w_2}}(w_1d-1)B_{j-k,\chi, \xi^{w_2}}^{(m-1)}(w_{1}y)
\\
& &=\underset{j=0}{\overset{n}{\sum}}\binom{n}{j}w_1^{j}
w_2^{n-j-1}B_{n-j, \chi,
\xi^{w_2}}^{(m)}(w_1x)\underset{k=0}{\overset{j}{\sum}}\binom{j}{k}T_{k,
\chi, \xi^{w_1} }(w_2d-1)B_{j-k, \chi,
\xi^{w_1}}^{(m-1)}(w_2y).\notag
\end{eqnarray*}
\end{theorem}

\begin{remark}
Taking $m=1$ and $y=0$ in $(2.7)$ derives the following identity :

\begin{eqnarray}
& &\underset{j=0}{\overset{n}{\sum}}\binom{n}{j}w_2^{j}
w_1^{n-j-1}B_{n-j, \chi, \xi^{w_1}}(w_2x)T_{j,
\chi, \xi^{w_2}}(w_1d-1)\\
& &\quad =\underset{j=0}{\overset{n}{\sum}}\binom{n}{j}w_1^{j}
w_2^{n-j-1}B_{n-j, \chi, \xi^{w_2}}(w_1x)T_{j, \chi, \xi^{w_1}
}(w_2d-1).\notag
\end{eqnarray}
\end{remark}

Moreover, if we take $x=0$ and $y=0$ in Theorem 1, then we have the
following identity for the generalized twisted Bernoulli numbers of
higher order.

\begin{corollary}
Let $d, w_1, w_2 \in \mathbb{N}$. For $n \in \mathbb{Z}_+$ and $m
\in \mathbb{N}$, we have
\begin{eqnarray*}
& &\underset{j=0}{\overset{n}{\sum}}\binom{n}{j}w_2^{j}
w_1^{n-j-1}B_{n-j, \, \chi,
\xi^{w_1}}^{(m)}\underset{k=0}{\overset{j}{\sum}}\binom{j}{k}T_{k,
\chi, \xi^{w_2}}(w_1d-1)B_{j-k,\chi, \xi^{w_2}}^{(m-1)}
\\
& &=\underset{j=0}{\overset{n}{\sum}}\binom{n}{j}w_1^{j}
w_2^{n-j-1}B_{n-j, \chi,
\xi^{w_2}}^{(m)}\underset{k=0}{\overset{j}{\sum}}\binom{j}{k}T_{k,
\chi, \xi^{w_1} }(w_2d-1)B_{j-k, \chi, \xi^{w_1}}^{(m-1)}.\notag
\end{eqnarray*}
\end{corollary}
We also note that taking $m=1$ in Corollary 2 shows the following
identity :
\begin{eqnarray}
& &\underset{j=0}{\overset{n}{\sum}}\binom{n}{j}w_2^{j}
w_1^{n-j-1}B_{n-j, \chi, \xi^{w_1}}T_{j,
\chi, \xi^{w_2}}(w_1d-1)\\
& &\quad =\underset{j=0}{\overset{n}{\sum}}\binom{n}{j}w_1^{j}
w_2^{n-j-1}B_{n-j, \chi, \xi^{w_2}}T_{j, \chi, \xi^{w_1}
}(w_2d-1).\notag
\end{eqnarray}

Now we will derive another interesting identities for the
generalized twisted Bernoulli numbers and polynomials of higher
order.  From $(\ref{1.15})$, $(\ref{2.2})$ and $(\ref{2.3})$, we can
derive that

\begin{eqnarray}
& &Y(m, \chi, \xi | w_1, w_2) \notag \\
&&=\frac{1}{w_1} \left(\sum_{i=0}^{w_1
d-1}\chi(i)\,\xi^{w_2i}\int_{X^m}(\underset{i=1}{\overset{m}{\prod}}\chi(x_i))\xi^{(\underset{i=1}{\overset{m}{\sum}}
x_i)w_1
}e^{(\underset{i=1}{\overset{m}{\sum}}x_i+\frac{w_2}{w_1}i+w_2
x)w_1t}dx_1 \cdots
 dx_m\right) \label{2.9}\\
& & \quad \times
 \left(\int_{X^{m-1}}(\underset{i=1}{\overset{m-1}{\prod}}\chi(x_i))\xi^{(\underset{i=1}{\overset{m-1}{\sum}}
x_i)w_2 }e^{(\underset{i=1}{\overset{m-1}{\sum}}x_i+w_1 y)w_2t}dx_1
 \cdots  dx_{m-1}\right) \notag\\
&&=\underset{n=0}{\overset{\infty}{\sum}}\left(\underset{k=0}{\overset{n}{\sum}}\binom{n}{k}w_1^{k-1}
w_2^{n-k}B_{n-k, \chi,
\xi^{w_2}}^{(m-1)}(w_{1}y)\underset{i=0}{\overset{w_1d
-1}{\sum}}\chi(i)\xi^{w_2 i}B_{k, \chi,
\xi^{w_1}}^{(m)}(w_{2}x+\frac{w_2}{w_1}i)\right)\frac{t^n}{n!}.\quad
\notag
\end{eqnarray}
From the symmetry property of $Y(m, \chi, \xi | w_1, w_2)$ in $w_1$
and $w_2$, we see that
\begin{eqnarray}
& &Y(m, \chi, \xi | w_1, w_2)  \label{2.10}\\\
&&=\underset{n=0}{\overset{\infty}{\sum}}\left(\underset{k=0}{\overset{n}{\sum}}\binom{n}{k}w_2^{k-1}
w_1^{n-k}B_{n-k,\chi, \xi^{w_1
}}^{(m-1)}(w_{2}y)\underset{i=0}{\overset{w_2d
-1}{\sum}}\chi(i)\xi^{w_1 i}B_{k, \chi, \xi^{w_2
}}^{(m)}(w_{1}x+\frac{w_1}{w_2}i)\right)\frac{t^n}{n!}. \quad \notag
\end{eqnarray}

Comparing the coefficients on the both sides of $(\ref{2.9})$ and
$(\ref{2.10})$, we obtain the following theorem which shows the
relationship between the power sums and the generalized twisted
Bernoulli polynomials.
\begin{theorem}
Let $d, w_1, w_2 \in \mathbb{N}$. For $n \in \mathbb{Z}_+$ and $m
\in \mathbb{N}$, we have
\begin{eqnarray*}
& &\underset{k=0}{\overset{n}{\sum}}\binom{n}{k}w_1^{k-1}
w_2^{n-k}B_{n-k,\chi,
\xi^{w_2}}^{(m-1)}(w_{1}y)\underset{i=0}{\overset{w_1d
-1}{\sum}}\chi(i)\xi^{w_2 i}B_{k,\chi, \xi^{w_1}}^{(m)}(w_{2}x+\frac{w_2}{w_1}i)  \\
& &=\underset{k=0}{\overset{n}{\sum}}\binom{n}{k}w_2^{k-1}
w_1^{n-k}B_{n-k, \chi,
\xi^{w_1}}^{(m-1)}(w_{2}y)\underset{i=0}{\overset{w_2d
-1}{\sum}}\chi(i)\xi^{w_1 i}B_{k, \chi,
\xi^{w_2}}^{(m)}(w_{1}x+\frac{w_1}{w_2}i).\notag
\end{eqnarray*}
\end{theorem}

\begin{remark}
Let $m=1$ and $y=0$ in Theorem 3. Then it follows that
\begin{eqnarray}
w_1^{n-1} \underset{i=0}{\overset{w_1d -1}{\sum}}\chi(i)B_{n,\chi,
\xi^{w_1}}(w_{2}x+\frac{w_2}{w_1}i) =w_2^{n-1}
\underset{i=0}{\overset{w_2d -1}{\sum}}\chi(i)B_{n,\chi,
\xi^{w_2}}(w_{1}x+\frac{w_1}{w_2}i).\label{2.11}
\end{eqnarray}
\end{remark}

Moreover, if we take $x=0$ and $y=0$ in Theorem 3, then we have the
following identity for the generalized twisted Bernoulli numbers of
higher order.

\begin{corollary}
Let $d, w_1, w_2 \in \mathbb{N}$. For $n \in \mathbb{Z}_+$ and $m
\in \mathbb{N}$, we have
\begin{eqnarray*}
& &\underset{k=0}{\overset{n}{\sum}}\binom{n}{k}w_1^{k-1}
w_2^{n-k}B_{n-k,\chi, \xi^{w_2}}^{(m-1)}\underset{i=0}{\overset{dw_1
-1}{\sum}}\chi(i)\xi^{w_2 i}B_{k,\chi, \xi^{w_1}}^{(m)}(\frac{w_2}{w_1}i) \\
& &=\underset{k=0}{\overset{n}{\sum}}\binom{n}{k}w_2^{k-1}
w_1^{n-k}B_{n-k, \chi,
\xi^{w_1}}^{(m-1)}\underset{i=0}{\overset{dw_2
-1}{\sum}}\chi(i)\xi^{w_1 i}B_{k, \chi,
\xi^{w_2}}^{(m)}(\frac{w_1}{w_2}i).\notag
\end{eqnarray*}
\end{corollary}

If we take $m=1$ in Corollary 4, we derive the identity for the
generalized twisted Bernoulli numbers : for $d, w_1, w_2 \in
\mathbb{N}$ and $n \in \mathbb{Z}_+$,
\begin{eqnarray}
w_1^{n-1}\sum_{i=0}^{dw_1-1}\chi(i)\xi^{w_2 i} B_{n, \chi,
\xi^{w_1}}( \frac{w_2}{w_1}i)
=w_2^{n-1}\sum_{i=0}^{dw_2-1}\chi(i)\xi^{w_1 i} B_{n, \chi,
\xi^{w_2}}(\frac{w_1}{w_2}i). \label{2.12}
\end{eqnarray}

\noindent \textbf{Acknowledgement.} This research was supported by
Kyungpook National University research fund 2008.

\end{document}